\documentclass[10 pt,leqno]{amsart}

\usepackage {amsfonts}
\usepackage{amsthm}
\usepackage{amssymb}
\usepackage{latexsym}
\usepackage{amsmath}
\usepackage{mathrsfs}
\usepackage{pifont}

\theoremstyle{plain}
\newtheorem{tw}{Theorem}[section]

\newtheorem {lem} [tw]{Lemma}

\newtheorem{cor}[tw]{Corollary}

\theoremstyle{definition}
\newtheorem {deft}[tw] {Definition}
\newtheorem {rem} [tw]{Remark}
\newtheorem {example} [tw]{Example}

\newcommand{\bc} {\Bbb C}

\newcommand{\bt} {\underline{k}}
\newcommand{\bm} {\underline{m}}
\newcommand{\bn}{\Bbb N}
\newcommand{\br}{\Bbb R}

\newcommand{\alg}{\mathcal{A}}
\newcommand{\hil}{\mathcal{H}}

\begin{document}
\author{Adam G. Skalski}
%\footnote{\emph{Permanent address of the author}.
\title{ On a classical scheme in noncommutative multiparameter ergodic theory}
\address{Department of Mathematics, University of \L\'{o}d\'{z}, ul.
Banacha 22, 90-238 \L\'{o}d\'{z}, Poland.}
\email{adskal@math.uni.lodz.pl}


\begin{abstract}  In the first part of the paper we describe
         the natural scheme for proving noncommutative individual ergodic
         theorems, generalize it for multiple sequences of measurable
         operators affiliated with a semifinite von Neumann algebra $M$, and
         apply it to theorems concerning unrestricted
         convergence of multiaverages.
         In the second part we prove convergence of ergodic averages
         induced by several maps satisfying specific recurrence relations,
         including so-called Multi Free Group Partial Sums.
         This is the multiindexed version of results obtained earlier jointly
         with V.I.Chilin and S.Litvinov.

\end{abstract}

\maketitle

Many interesting cases considered in ergodic theory,
both classical and quantum, can be cast in the common framework. One usually
deals with a family of transformations (representing or evolution of a system either
averages of some quantities over periods of time), and asks questions about
the convergence of these maps. Positive answers to these questions can be understood as the existence
of some limit behaviour of a system, or of a mean value of a given quantity.
The domain and range of these maps, the types of convergence, and the sense attributed
to them, all depend on the specific problem. However, often one can assume that the
evolution/averaging is linear, and that the domain of definition of our maps
is some normed space $B$ - clasically this might be the space of
integrable functions over a probability space; in the quantum context it might be
a $C^*$-algebra, a von Neumann algebra or a noncommutative $L^p$-space.
In this paper we will be especially interested in the multiparameter case,
corresponding physically to the existence of several (not necessarily independent)
evolutions of our system. We shall work in discrete time, and investigate
behaviour at infinity.

    The aim of this paper is to present applications of the well-known
classical scheme of proving individual ergodic theorems in the
noncommutative context. After establishing some necessary notations in the
 introductory section, in Section 1 we describe how to extend
this scheme to multisequences of maps acting in von Neumann algebras,
as was done in \cite{Ban}, \cite{subseq} and \cite{wspolna} for sequences
indexed by one parameter. Section 2 collects and briefly summarizes
known facts concerning unrestricted convergence of multiaverages and shows how to
reprove them using the aforementioned scheme. Finally in Section 3
we present a few ergodic theorems on averages induced by  several families
of maps satisfying specific recurrence relations (of which the so-called Multi
Free Group Actions are special examples). This is a multiparameter extension
of results established using similar methods in \cite{wspolna}, and derives
from earlier work of A.Nevo, E.Stein and T.Walker.

% \setcounter {section}{-1}
% \section{ Basic notions and definitions}

%\hspace{0.7 cm}
Let $d$ be a fixed positive integer.
All multiindices will be underlined and will usually  belong to ${\bn_0}^d$
or $\bn^d$,
where $\bn_0 = \bn \cup \{0\}$. When $\bt = \{ k_1, \ldots, k_d \}$
we will write $ \min {\bt} = \min \{k_1, \ldots, k_d \}$,
$\max {\bt} = \max \{k_1, \ldots, k_d \}.$
Below we recall the notion of unrestricted convergence (convergence in Pringsheim's sense)
for  a multisequence.

\begin{deft}
\label{max} Let $(\xi_{\bt})_{\bt \in{N_0} ^d}$ be a multisequence of real numbers. We say that it converges to
$\eta \in \br$ in Pringsheim's sense when for each $\epsilon >0$ there exists $n \in \bn$
such that  $ | \xi_{\bt} - \eta |< \epsilon $ whenever $\bt \in\bn_0^d$, $\min{\bt} \geq n$.
\end{deft}

Let $M$ be a semifinite von Neumann algebra with a faithful normal semifinite trace
$\tau$ (in some places we will loosen or strengthen the assumptions on $M$ and $\tau$).
Its positive part will be written as $M^+$, its hermitian part as $M_{\rm sa}$
and the lattice of all
projections belonging to $M$ as $P_M$. We will also denote $1-p$ by $p^{\perp}$ for $p \in P_ M$.
By $\widetilde{M}$ we shall denote the space of all measurable operators
affiliated with $M$, by $L^1(M)$ and $L^2(M)$ respectively the spaces of integrable
and square-integrable operators (see \cite{Nelson}). The space $\widetilde{M}$ will be equipped
with the topology of convergence in measure. Apart from the standard
convergences (in norm or in measure) one can define in all these spaces various equivalents
of the classical almost everywhere convergence. In our paper we will basically use two of
them, the almost uniform convergence introduced by E.Lance in \cite{Lance} and
the bilateral almost uniform convergence introduced by F.J.Yeadon in
\cite{Yeadon}:

\begin{deft}         \label{auconv}
A sequence $(x_n)_{n=1}^{\infty}$ of operators of $\widetilde{M}$
is almost uniformly (a.u.) convergent to $x \in \widetilde{M}$ if for each $\epsilon >0 $
there exists $p \in P_M$, $ \tau (p^{\perp}) < \epsilon$ and such that
\[ \left\| (x_n -x) p \right\|_{\infty} \stackrel{n \longrightarrow  \infty}
  {\longrightarrow} 0.\]

A sequence $(x_n)_{n=1}^{\infty}$ of operators of $\widetilde{M}$
is bilaterally almost uniformly  (b.a.u.) convergent to $x \in \widetilde{M}$ if for each
$\epsilon >0 $ there exists $p \in P_M$, $ \tau (p^{\perp}) < \epsilon$ and  such that
\[ \left\| p (x_n -x) p \right\|_{\infty} \stackrel{n \longrightarrow  \infty}
  {\longrightarrow} 0.\]
\end{deft}

%\vspace{0.3 cm}

We mention below several other counterparts of classical properties holding a.e..
Let $B$ be a linear space.

\begin{deft} \label{subadd} A map $a: B \longrightarrow \widetilde{M}$ is
bilaterally almost
 uniformly subadditive (b.a.u.\ subadditive) if
for all $\epsilon >0 $ there exists $p \in P_M$, $ \tau (p^{\perp}) < \epsilon$ and such that
\[ \left\|p a (x+y) p \right\|_{\infty} \leq \left\|p a (x) p \right\|_{\infty}
+ \left\|p a (y) p \right\|_{\infty}. \]

A map $a: B \longrightarrow \widetilde{M}$ is bilaterally almost
 uniformly homogenous (b.a.u.\ homogenous) if
for all  $\epsilon >0 $ there exists $p \in P_M$, $ \tau (p^{\perp}) < \epsilon$ and such that
 for all $x \in B$, $\alpha \in \br$
\[ \left\|p a (\alpha x) p \right\|_{\infty} = | \alpha | \left\|p a (x) p \right\|_{\infty}. \]
\end{deft}

%\vspace{0.6 cm}
The basic notion used in the following is that of \it kernel \rm (also called \it absolute
 contraction \rm).

\begin{deft}
A linear map $\alpha : L^1(M) \longrightarrow L^1(M)$ is a kernel if
it is a positive contraction:
\[ \forall_ {x\in L^1(M) \cap M} { \; \; 0 \leq x \leq I \Longrightarrow 0 \leq \alpha(X) \leq I}\]
and has the property
\[ \forall_ {x\in L^1(M)^, \, x\geq 0} {\; \; \tau (\alpha(x)) \leq \tau (x) }.\]
\end{deft}

Each kernel defines uniquely a $\sigma$-weakly continuous map transforming $M$
into $M$ coinciding with the original map on $M \cap L^1(M)$ and usually
denoted by the same letter. Moreover, for each kernel $\alpha$ there
exists a (unique) \it adjoint kernel \rm $\alpha^{\star}$ such that for all
$a \in L^1(M)$, $b \in M$
\[ \tau (\alpha(a) b ) = \tau (a \alpha^{\star}(b)) \]
(the proof of all these facts can be found in \cite{Yeadon}).

%\vspace{0.5 cm}

We will also need a specific type of dense subsets of a Banach space.

\begin{deft}
Let $(B, \ \| \cdot \|, \ \geq )$ be an ordered real Banach space with the closed
convex cone $B_+$, $B=B_+-B_+$. A subset $B_0 \subset B_+$ is said to be
minorantly dense in $B_+$ if for every $x\in B_+$ there is a sequence
$\{ x_n \}$ in $B_0$ such that $x_n \leq x$ for each $n$, and
$\| x-x_n \| \to 0$ as $n\to \infty$. For example, $M^+ \cap L^1(M) \cap L^2(M)$
is a minorantly dense subset of $L^1(M)_{\rm sa}$.
\end{deft}

%\vspace{1cm}

\section{ A classical scheme for proving individual ergodic theorems in the noncommutative
 context}

%\hspace{0.7 cm}
 According to the description in the introduction we consider the following
situation: let $B$ be a Banach space, and let $(a_{\bt})_{\bt \in \bn_0^d}$ be
a family of linear maps acting from $B$ to the algebra $\widetilde{M}$ of all
measurable operators affiliated with some semifinite von Neumann algebra $M$ (in the
classical case $\widetilde{M}$ is replaced by the algebra $L_0(X,\mu)$ of all
measurable functions on a measure space $(X, \mu)$). The typical way of
proving pointwise (or b.a.u.) convergence of a family $(a_{\bt} (x))_{\bt \in \bn_0^d}$
for each $x \in B$ consists of three parts:

\begin{itemize}
\item [\bf{A}] prove that we have the desired convergence for each $x \in B_0$
- some dense subset of $B$
\item [\bf{B}] establish some estimate of the maximal type (separately for each $x \in B$)
\item [\bf{C}] deduce the convergence for all $x\in B$ , using {\bf A}, {\bf B}
 and  continuity properties of the maps $a_{\bt}$
\end{itemize}

We will now briefly describe each of these steps.

In most cases  (and commutative and noncommutative) $B$ is the $L^1$-space, and $B_0$
is either the $L^2$-space or the span of indicator functions/projections.
The technique usually used in obtaining {\bf A} is
first to prove some kind of  mean ergodic theorem (in the spirit of von
Neumann's theorem on averages of contractions in a Hilbert space) and then use
it to deduce the pointwise convergence - see section 3.

As concerns {\bf B}, basically all of the known maximal lemmas used in the
noncommutative ergodic theory are based on the maximal lemma of F.Yeadon
(\cite{Yeadon}).
For our purposes it is convenient to formulate the version of this lemma
 proved in \cite{Petz1}, for a sequence of operators. In fact it was proved in a
greater generality which we shall use in Remark \ref{comment}.

\begin{tw}
If $\alpha: M \longrightarrow M$ is a kernel, $(x_m)_{m=1}^{\infty}$ is a sequence of operators
belonging to $L^1(M)^+$ and $(\epsilon_m)_{m=1}^{\infty}$ is a sequence of positive real numbers
(estimation numbers),
then there exists a projection $p \in P_M$ such that
\[\tau(p^\perp) \leq 2 \sum_{m=1}^{\infty}{\epsilon_m^{-1} \tau(x_m)},\]
\[ \|p \left(\frac{1}{r}\sum_{k=0}^{r-1} {\alpha^k(x_m)} \right) p \|_{\infty} \leq 2 \epsilon_m, \; \; r,m \in \bn.\]
\end{tw}

For the case of finite $M$, and operators $x_m\in M^+\cap L^1(M)$, one can replace the
above estimates with the one-sided version (using the Kadison
inequality).

As we are interested here in the case of several kernels, we will need the
following extension of the above maximal lemma:

\begin{tw}  Let $ \alpha_1,\ldots, \alpha_d:M \to M$ be mutually commuting
kernels.
 Then there exists a constant $\chi_d >0$ such that for every sequence of operators
$(x_m)_{m=1}^{\infty}$ belonging to $L^1(M)^+$ and
$(\epsilon_m)_{m=1}^{\infty}$ - a sequence of positive real numbers
(estimation numbers), there exists $p\in P(M)$ such that
\[\tau(p^\perp) \leq 2 \sum_{m=1}^{\infty}{\epsilon_m^{-1} \tau(x_m)},\]
\[ \|p \left(\frac{1}{n^d} \sum_{k_1=0}^{n-1}\ldots \sum_{k_d=0}^{n-1}
{\alpha_1^{k_1} \circ \ldots \circ \alpha_d^{k_d}(x_m)} \right) p \|_{\infty}
 \leq 2 \chi_d \epsilon_m, \; \; n,m \in \bn.\]
\label{Petzmax}
\end{tw}

\begin{proof} The proof of the theorem stems from a fact proved by A.Brunel in \cite{Brunel} for
commuting contractions acting in the classical $L^1$. Repeating his reasoning we
can find for each $\bt \in \bn_0^d$ a nonnegative number $a(\bt)$ such that

(i) $ \sum_{\bt \in \bn_0^d} {a(\bt)} =1,$

(ii) the mapping $U$ defined by $U = \sum_{\bt \in \bn_0^d} {a(\bt)  \alpha_1^{k_1}
  \circ \ldots \circ \alpha_d^{k_d}}$ satisfies the following inequality:

\[\frac{1}{n^d} \sum_{k_1=0}^{n-1}\ldots \sum_{k_d=0}^{n-1}
{\alpha_1^{k_1} \circ \ldots \circ \alpha_d^{k_d} (x) }\leq
 \frac{\chi_d}{n_d} \sum_{j=0}^{n_d-1} {U^j (x)},\]
for any $x\in L^1 (M)$, $n \in \bn$, where $\chi_d >0$ depends only on $d$ and
$n_d \in \bn$  depends only on $d$ and $n$.

As this is clear that $U$ is also a kernel, the above version of
Yeadon's theorem ends the proof.
\end{proof}

\begin{rem}
The above theorem remains true if $M$ is any von Neumann algebra with a n.s.f.\
weight $\phi$, $\alpha_i:i=1,\ldots d$ are positive linear maps acting
in $M$ such that for each $x \in M, 0 \leq x
\leq I,$ we have $\alpha_i(x) \leq I$, for each $x \in M^+ $ we have
$\phi (\alpha_i(x)) \leq \phi(x)$,          \label{comment}
and the sequence $(x_m)_{m=1}^{\infty}$ consists of operators in $M^+$.
Naturally, we have to replace everywhere $\tau$ by $\phi$.
\end{rem}

The classical tool for {\bf C} is the Banach Principle, established already in
1926. Its noncommutative
generalization was proved by M.Goldstein and S.Litvinov in \cite{Ban}
for the quasi uniform convergence, and then also by V. Chilin, S. Litvinov
and the author (\cite{wspolna}) for the b.a.u.\ convergence.
Here we need the extension of this result for multisequences (due in the
classical case to F.Moricz \cite{Moricz}).  Because of some technical
subtleties we need to work with minorantly dense subsets of a Banach space.

\begin{tw} \label{Banach}
        Let $d\in \bn$, let $B$ be an ordered real Banach space with the
closed convex cone $B_+$, $B_+ - B_+ = B$,
 and for each $\bt : =(k_1,\ldots , k_d) \in {\bn_0} ^d$
let $a_{\bt} : B \longrightarrow \widetilde{M}$ be a continuous positive
linear map. Assume that the following
conditions are satisfied: \\
\rm (i) \it for each $b \in B_+$ and $\delta >0$ there exists $y \in M^+$,
$0\neq y \leq I$ and $K \in \bn$ such that
\[ {\sup_{ \rm min \it {\underline{k} \geq K}} \left\| y a_{\bt} (b) y \right\|_{\infty}} < \infty \]
\hspace {0.5cm} and $\tau(I-y) \leq \delta,$ \\
\rm (ii) \it there exists $B_0$, a minorantly dense subset of $B_+$ such
that for each $b\in B_0$ the operators
$a_{\bt} (b) - a_{\underline{m}} (b)$
b.a.u. converge to 0 as $\bt, \underline{m} \longrightarrow \infty$, in Pringsheim's sense. \\
\hspace{0.8 cm} Then for each $b \in B$, $a_{\bt} (b) $ is b.a.u. convergent to some element of $\widetilde{M}$
as $\bt \longrightarrow \infty$ in Pringsheim's sense.
\end{tw}

\begin{proof} The method of proof is typical; it is based on the Baire Category Theorem, and is
reminiscent of that adopted in \cite{Ban} and \cite{wspolna}. Whenever we say that some multisequence
of real numbers converges to a real number, it is to be understood in Pringsheim's sense.

Let's fix $\epsilon > 0$. For each $j, L, K  \in \bn$ we put
\[ \epsilon_j = \epsilon / 2^{j+3},\]
\begin{equation} {\large B_{j,L,K} = \left\{ b \in B_+ : \, \exists_{y \in M^+, \, y \leq 1, \, \tau (1-y) \leq \epsilon_j}
 \, \forall_{ \bt \in {N_0}^d, \,  \min{\bt} \geq K} \, \left\|a_{\bt} (b)^{\frac{1}{2}} y \right\|_{\infty} \leq L \right\}.
 \label{def1} } \end{equation}
Fix now $j \in \bn$. It is easy to see that
\[ B_+ = \bigcup_{L,K =1}^{\infty} {B_{j,L,K}}. \]
Moreover, using $\sigma$-weak compactness of the unit ball of $M$ one can
show (exactly as was done in \cite{wspolna} for the sets $X_{L,k}$) that each of the
sets $B_{j,L,K}$ is closed. Once this has been done, the Baire Theorem allows us to infer that there
exist $L_j, K_j \in \bn$, $b_j \in B_+$ and $\delta_j >0$ such that $B_{j, L_j, K_j}$ contains
the ball with centre in $b_j$ and radius $\delta_j$. This means that for any
$b \in B_+$, $\| b- b_j \| \leq \delta_j$ there exists $y_{b,j} \in M^+$ satisfying conditions
mentioned in (\ref{def1}). Let
\[ y_{b,j} = \int_{0}^{1} { \lambda d E_{b,j}(\lambda)} \]
be the spectral decomposition of $y_{b,j}$. Define
\[ g_{b,j} = 1 - E_{b,j} (\frac{1}{2}).\]
We have
\begin{equation} \tau( g_{b,j}^{\perp}) \leq 2 \epsilon_j, \label{proj1}\end{equation}
\begin{equation} \|  g_{b,j} a_{\bt} (b) g_{b,j} \|_{\infty} \leq 2 L_j^2  \label{proj2}\end{equation}
for all $\bt \in {\bn_o}^d$, $ \min {\bt} \geq K_j$.
Whenever $c \in B_+$, $ \| c \| \leq  \frac{\delta_j}{ 4 j L_j^2}$, putting
\begin{equation} f_{c,j} = g_{ b_j , j } \wedge g_{ b_j - 4 j L_j^2 c, j }  \label{proj3}\end{equation}
and using (\ref{proj1}), (\ref{proj2}) we get for all $\bt \in {\bn_o}^d$,
$ \min {\bt} \geq K_j$
\[ \| f_{c,j} a_{\bt} (c) f_{c,j} \|_{\infty} = (4 L_j^2 j)^{-1}  \| f_{c,j} a_{\bt} (4 L_j^2 j c) f_{c,j} \|_{\infty} \leq  \]
\[ \leq   (4 L_j^2 j)^{-1} \left( \| f_{c,j} a_{\bt} ( b_j) f_{c,j} \|_{\infty} + \| f_{c,j} a_{\bt} ( b_j - 4 j L_j^2 c, j) f_{c,j}
   \|_{\infty} \right) \leq \]
\begin{equation} \leq (4 L_j^2 j)^{-1} \left( \| g_{b_j,j} a_{\bt} ( b_j) g_{b_j,j} \|_{\infty} + \|  g_{b_j - 4 j L_j^2 c,j} a_{\bt}
( b_j - 4 j L_j^2 c, j) g_{b_j - 4 j L_j^2 c,j}
   \|_{\infty} \right) \leq  \frac{1}{j} . \label{ proj4} \end{equation}

Now let $b$ be a fixed element of $B_+$. There exists a sequence $(c_j)_{j=1}^{\infty}$ of elements of $B$
such that for each $j \in \bn$
\[ \left\| c_j \right\| \leq \frac{\delta_j}{ 4 j L_j^2 }, \; \, b+c_j \in B_0 .\]
We choose a sequence $(p_j)_{j=1}^{\infty}$ of projections from $M$ such that:
for each $ j \in \bn$
\[  \tau (p_{j}^{\perp}) < \epsilon_{j} , \; \left\| (a_{\bt} p_j (b+ c_{j}) - a_{\bm} (b+ c_{j})) p_{j} \right\|_{\infty}
 \stackrel {{\bt}, {\bm} \rightarrow \infty }{\longrightarrow } 0. \]
Put
\[ q = \bigwedge_{j=1}^{\infty} {p_j} \wedge \bigwedge_{j=1}^{\infty} {f_{b+c_j, j}}. \]
The definition of $(p_j)_{j=1}^{\infty}$, together with (\ref{proj1}) and  (\ref{proj3}), gives
\[ \tau (q^{\perp}) \leq \sum_{j=1}^{\infty}{\epsilon_j} + \sum_{j=1}^{\infty}{4 \epsilon_j} < \epsilon. \]

It remains to prove that $\| q(a_{\bt} (b)- a_{\bm}(b)) q \|_\infty $ tends to zero as
$\bt, \bm \longrightarrow \infty$. Fix $\delta >0$, and let $j \in \bn$ be such that
$ \delta > 3 j^{-1}$. Then
\begin{align*} \| q & (a_{\bt}(b) -   a_{\bm}(b)) q \|_\infty  \leq \\
&\| q (a_{\bt}(b+c_j) - a_{\bm}(b+c_j)) q \|_\infty
  + \| q a_{\bt}(c_j) q \|_\infty + \|q a_{\bm}(c_j) q \|_\infty \leq \\
 &\| p_j (a_{\bt}(b+c_j) - a_{\bm}(b+c_j))  p_j \|_\infty
  + \| f_{c_j, j} a_{\bt}(c_j) f_{c_j, j} \|_\infty +
  \| f_{c_j, j} a_{\bm}(c_j) f_{c_j, j} \|_\infty < \delta \end{align*}
for all ${\bt}, {\bm} \in {\bn_0}^d$ such that $\min{\bt} \geq K$ and $\min{\bm} \geq K$,
moreover $K$
depends only on $b$ and $j$ (so actually only on $b$ and $\delta$). As the algebra
$\widetilde{M}$ is complete with respect to the topology of b.a.u. convergence (Theorem 2.3 of
\cite{wspolna}),
this ends the proof of the desired convergence for any $b \in B_+$. The general case follows
 immediately. \end{proof}

%\vspace{0.3 cm}
The above theorem remains true when one replaces throughout (both in the assumptions
and in the hypothesis) convergence
in Pringsheim's sense by the so-called 'maximal' convergence.  Obviously one also has
to reformulate properly the condition on the b.a.u.\ boundedness of the maps considered.
Moreover, when one replaces b.a.u.\ convergence by quasi uniform convergence,
one may prove the theorem considering a dense subset of a Banach space
(instead of a minorantly dense subset of an ordered Banach space).
A careful reader would also notice that actually it is enough to assume that each map
$a_{\bt} : B\longrightarrow \widetilde{M}$ is positive, b.a.u.\ homogeneous and subadditive (see definitions
(\ref{subadd})) instead of assuming that the $a_{\bt}$'s are linear.

   In our context the technical complication (using additionally the order structure
in a given Banach space) will not be a serious obstacle.
Working with complex Banach spaces of operators (say $L^1(M)$) we can first concentrate
on the selfadjoint parts of them, and then use the existence of the convenient
decomposition of a given operator into its real and imaginary part to
conclude the convergence of an investigated sequence. This kind of reasoning
will be further used without any comments.

%\vspace{1cm}

\section{ Unrestricted convergence of  multiparameter averages}

%\hspace{0.7 cm}
The theorem below appeared first in \cite{GoGr} and then was
also mentioned in \cite{Junge}. Here we show how to deduce it
immediately from the maximal lemma established in the
second-mentioned paper and the older result of D.Petz.

\begin{tw} [\cite{GoGr}, \cite{Junge}]
Let $d \in \bn$, $p \in (1, \infty)$, $\alpha_i : L^1(M) \longrightarrow L^1(M)$ ($i=1,\ldots, d$)
 be kernels. For each $y \in L^p(M)$ denote by $\Phi_i (y)$ the norm limit
of the sequence $\left(\frac{1}{n} \sum_{k=0}^{n-1} {\alpha_i^k (y)}
\right)_{n=1}^{\infty}$ (which exists by the reflexivity of $L^p$-space). Then  for each $x\in L^p(m)$
the multisequence  $\left(s_{\bt} (x)\right)_{\bt \in \bn^d}$, where
\[s_{\bt}(x) =  \frac{1}{k_1 \ldots k_d} \sum_{i_1 =0}^{k_1-1} \ldots  \sum_{i_d =0}^{k_d-1}
   {\alpha_1^{i_1} \circ \ldots \circ \alpha_d^{k_d} (x)},\]
 b.a.u.\ converges to the operator
$\Phi_1(\ldots(\Phi_d(x))\ldots)$ as $\bt\longrightarrow \infty$ in Pringsheim's sense.
 \label{pointerg} \end{tw}

\begin{proof}
       In \cite{Junge} it was proved that for each $y \in L^p(M)^+$,
and each kernel $\beta:M \longrightarrow M$ there exists an operator
$\tilde{y}$ such that for all $n\in \bn$
\begin{equation} \frac{1}{n}\sum_{i=0}^{n-1} {\beta^i (y)} \leq \tilde{y}. \label{estim}\end{equation}
In our context this immediately implies that there exists $\tilde{x}$
such that
\[s_{\bt}(x) \leq \tilde{x}. \]
This is sufficient for part {\bf B} of our scheme.
Putting $B_0 = M^+ \cap L^1(M) \cap L^2(M)$ and using theorem
4 of \cite{Petz1} we obtain part {\bf A}. Theorem \ref{Banach} shows that the
multisequence considered is b.a.u.\ convergent, and standard reasoning allows us to
conclude that the b.a.u.\ limit is equal to the norm limit.
\end{proof}

\begin{rem}
In \cite{GoGr} a formula analogous to formula
(\ref{estim}) was obtained for some
$\tilde{y} \in L^{p-\epsilon} (M)$ (for any given sufficiently small $\epsilon>0$).
 This clearly also suffices to conclude the proof in the same way as was done above.
\end{rem}

We would also like to briefly describe the situation concerning
norm convergence. Here nothing depends on the number of
 kernels considered, the only important factor is the finiteness of the
trace. This is illustrated by the following basic example

\begin {example} \rm
    Let $M= L^{\infty} (\br)$, with trace given by the Lebesgue integral, and let $\alpha$ be the
standard shift operator,

$$\left( \alpha(f) \right) (t) = f (t-1)$$
for all $f \in L^1(\br)$, $ t \in \br$. It is clear that $\alpha$ is a kernel, and
$$ s_k (\chi_{(0,1)} )  = \frac{1}{k}
\sum_{j=0}^{k-1} {\alpha^j(\chi_{(0,1)})} \stackrel{a.e.}{\longrightarrow} 0, $$
$$ \left\| s_k (\chi_{(0,1)} ) \right\|_{1}  = 1 $$
for all $k\in \bn$   (by $\chi_{(0,1)}$ we understand the characteristic function of the interval $(0,1)$).

\end{example}

The situation described above cannot happen when $M=L^{\infty} (X,\mu)$ and $(X,\mu)$ is a finite
measure space. On the algebraic level this corresponds to the finiteness of the trace on the
algebra $M$.

%\vspace{0.5cm}
\begin{tw}            \label{mean}
   Let $M$ be a von Neumann algebra with a faithful normal finite trace $\tau$.
\noindent Let $d \in \bn$, $\alpha_i : L^1(M) \longrightarrow L^1(M) $ ($i=1,\ldots, d$)
 be kernels and let
$x\in L^1(M)$. Then the multisequence  $s_{\bt} (x) $ converges in $L^1$-norm  as $\bt\longrightarrow \infty$ in the Pringsheim's sense.
\label{normerg} \end{tw}

\begin{proof}
The proof follows by induction
with respect to the number of kernels.  For each $r \in \{1,\ldots,d\}$ and $x \in L^1(M)$ if only the sequence $(s_{r}^n(x))_{n=1}^{\infty}$
is convergent in $L^1$-norm its limit will be denoted again
by $\Phi_r(x)$ (in the case of one kernel we shall write simply $\Phi(x)$).

Let $\alpha: L^1(M) \longrightarrow L^1(M)$ be a kernel. For each $n \in \bn$
we define the set
\[ A_n = \{ y\in L^1(M): \, y^{\star}=y, \, -n I \leq y \leq n I \}.\]
One can see that if $(p_k)_{k=1}^{\infty}$ is a sequence of mutually orthogonal projections
in $P_M$, and $y\in A_n$ then
\[ \left| \tau (y p_k) \right| \leq \left| \tau (n p_k) \right| \leq n
 \left| \tau (\sum_{i=k}^{\infty} {p_i} ) \right|,\]
so the expression on the left side of the above inequality tends to $0$, as $k$ tends to
$\infty$, uniformly with respect to $y$. The theorem II.2 of \cite{Akemann} allows us to infer that
$A_n$ is weakly relatively compact. As $A_n$ is a convex and norm-closed subset of $L^1 (M)$,
Mazur's theorem shows that it is actually weakly compact.
    Now we can apply Theorem 2.1.1 of \cite{Krengel} (notice that $\alpha: A_n \longrightarrow
 A_n$) to conclude that for each $x \in A_n$ there exists $\Phi(x) \in A_n$ such that
$s_k (x) \stackrel{k\longrightarrow \infty} {\longrightarrow} \Phi(x)$,
$\alpha(\Phi(x)) =\Phi(x)$. As the set  $ \bigcup_{n=1}^{\infty} A_n$ is norm-dense
in the hermitian part of $L^1(M)$, and each operator in $L^1(M)$ can be expressed
as a sum of two hermitian operators, the proof of the theorem for the case
$d=1$ is finished.

Assume now that we know that the theorem holds for $d-1$. Then we can write
the following inequalities:
\[ \| \Phi_1(\ldots(\Phi_d(x))\ldots) - s_{\bt}(x)  \|_1  \leq\]
\[ \leq
  \| \Phi_1(\ldots(\Phi_{d-1}(\Phi_d(x)))\ldots)  - s_{\bm}(\Phi_d(x)) \|_1
  + \| s_{\bm} (\Phi_d(x) - \alpha_{d}^{k_d}(x)) \|_1 \leq \]
\[ \leq \| \Phi_1(\ldots(\Phi_{d-1}(\Phi_d(x)))\ldots)  - s_{\bm}(\Phi_d(x)) \|_1
  + \|  (\Phi_d(x) - \alpha_{d}^{k_d}(x) \|_1, \]
where $\bm \in \bn_{0}^{d-1}$, $\bm = \{k_1, \ldots, k_{d-1}\}$, and we used
the fact that each kernel is a contraction. It is easily seen that the first
part of the above expression tends to 0, as $\bm$ tends to infinity
in Pringsheim's sense (by the induction assumption), and the same can be said
about the second part as $k_d \longrightarrow \infty$.
\end{proof}

\begin{rem}
     A version of the above theorem for one kernel is due to C.Radin (\cite{Radin}). However
since our assumptions are slightly different, we do not need to introduce the abstract
notion of a unit in a predual of a von Neumann algebra. Moreover, we give a more
detailed proof.

As was mentioned in Theorem \ref{pointerg}, for $p \in (1,\infty)$ the norm convergence of (multi)averages follows
immediately from the reflexivity of the space in question. The same is true
in the non-tracial situation.
\end{rem}

\section{ Ergodic theorems for multiparameter Free Group Actions}

For the sake of clarity we shall restrict ourselves to the case
of $d=2$ throughout this section - all the results hold for general $d\in \bn$.
We begin with a general fact concerning the strong convergence of
averages, formulated in the spirit of von Neumann's ergodic theorem.

 Let for any $p \in (0,1)$
\[D_p= \left\{ z \in \bc : | \sqrt{z+4p-4p^2} + \sqrt{z-4p-4p^2}| \leq 2 \sqrt{p}, \right. \]
               \[ \left. |\sqrt{z+4p-4p^2} - \sqrt{z-4p-4p^2} |\leq 2 \sqrt{p} \right\}. \]

The next theorem is an easy consequence of the following lemma:

\begin{lem}[\cite{Walker}] \label{comp}
          Assume that $p \in (0,1]$ and  $(f_n)_{n=1}^{\infty}$ is the sequence of functions
on the complex plane such that $f_0(z)=1$, $f_1 (z) = z $ and $z f_n(z) = p f_{n-1}(z) +
(1-p) f_{n+1} (z)$ for $n\geq 1$ ($z \in \bc$). Then
$\frac{1}{n} \sum_{k=0}^{n-1} f_k(z)$ converges pointwise iff $z \in D_p$. It converges to zero
on $D_p \setminus\{1\}$. Moreover there exists $N\in \bn$ such that for all $n \geq N$ and
$z \in D_p$  $|\frac{1}{n} \sum_{k=0}^{n-1} f_k(z)| \leq 1.$

\end{lem}

%\vspace{0.5 cm}

\begin{tw} Let $H$ be a Hilbert space, $x_{0,1}, x_{1,0}$ - commuting normal operators in $B(H)$
whose spectra
are respectively subsets of $D_{p_1}$ and $D_{p_2}$ for some.
$p_1, p_2 \in (0, 1]$. Let $x_{0,0} =I$ and $x_{m,n}$ (for $m+n\geq 2$),
 operators satisfying the relations $x_{1,0} x_{m,n} = p_1 x_{m+1,n} + (1-p_1) x_{m-1,n}$ be
 $x_{0,1} x_{m,n} = p_2 x_{m,n+1} + (1-p_2) x_{m,n-1}$. Then the multisequence
  $\left(\frac{1}{k_1 k_2} \sum_{m=0}^{k_1-1} \sum_{n=0}^{k_2-1}{x_{m,n}}\right)_{k \in \bn^2}$
   converges in Pringsheim's
  sense to the projection
$P$ onto the set $\{\eta \in H : x_1 \eta = x_2 \eta = \eta \}$.       \label{vNWalker}
\end{tw}

\begin{proof}
Let $E_{\lambda, \mu}$ denote the spectral measure on $D_{p_1} \times D_{p_2}$ corresponding to
the pair $x_{0,1}, x_{1,0}$. It is easy to see that
\[ \frac{1}{k_1 k_2} \sum_{m=0}^{k_1-1} \sum_{n=0}^{k_2-1}{x_{m,n}} =
\frac{1}{k_1 k_2} \sum_{m=0}^{k_1-1} \sum_{n=0}^{k_2-1} \int_{D_{p_1} \times D_{p_2}}
{ f_m^{(1)} (\lambda) f_n^{(2)} (\mu) dE_{\lambda, \mu}},\]
where by  $(f^{(1)}_n)_{n=1}^{\infty}$, $(f^{(2)}_n)_{n=1}^{\infty}$ we understand the sequences
introduced in lemma \ref{comp} with respectively $p=p_1$ and $p=p_2$. Therefore
\[ \frac{1}{k_1 k_2} \sum_{m=0}^{k_1-1} \sum_{n=0}^{k_2-1}{x_{m,n}} - P =
  \int_{D_{p_1} \times D_{p_2}\setminus \{(1,1)\}}
{ \frac{1}{k_1 k_2} \sum_{m=0}^{k_1-1} \sum_{n=0}^{k_2-1} f_m^{(1)} (\lambda) f_n^{(2)} (\mu) dE_{\lambda, \mu}},\]
and the desired strong convergence follows from standard properties of spectral integrals.
\end{proof}

The following notations will be used: let $\alg$ be a von Neumann algebra
with a faithful normal semifinite weight $\phi$. Further let
\[ \mathcal{N}_{\phi} = \{ A\in \alg : \phi (A^{\star} A) < \infty\}, \; \; \alg_0=
   \mathcal{N}_{\phi}^{\star} \cap  \mathcal{N}_{\phi},\]
and let $\hil_{\phi}$ be the Hilbert space completion of $\alg_0$  (with respect to the scalar
product $\langle A , B\rangle_{\phi} = \phi(B^{\star}A)$).  We will write $\Lambda_{\phi}$
for the canonical injection of $\alg_0$ in $\hil_{\phi}$, and $\pi_{\phi} : \alg \longrightarrow
B(\hil)$ for the faithful normal representation such that for all $A \in \alg, B \in \alg_0$
\[ \pi_{\phi} (A) (\Lambda_{\phi}(B)) = \Lambda_{\phi}(AB)\]
(left regular representation). We will also occasionally use the standard language of Hilbert algebras,
as in \cite{Takesaki}.

Let us describe the averages we will consider.
Fix real numbers $p_1,p_2 \in (0,1]$. Let $\sigma_{0,1}$, $\sigma_{1,0}$ be normal,
completely positive, unital, $\phi$-invariant and commuting maps acting on the algebra $\alg$.
Moreover let $\sigma_{m,n}$ (for $m.n \in \bn_0, m+n\geq 2$) be positive maps acting on $\alg$
defined recursively by the following relations:
\[\sigma_{1,0} \sigma_{m,n} = p_1 \sigma_{m+1,n} + (1 - p_1) \sigma_{m-1,n},\]
\[\sigma_{0,1} \sigma_{m,n} = p_2 \sigma_{m,n+1} + (1 - p_2) \sigma_{m,n-1},\]
 where
$\sigma_{0,0} = \rm Id_{\alg}$. Clearly         for all $m,n \in \bn_0$
\[ \sigma_{m,n} = \sigma_{0,n} \circ \sigma_{0,m}.\]
 We will write for each $k,k_1,k_2\in \bn$,
\[ S_{k_1,k_2} = \frac{1}{k_1 k_2}\sum_{m=0}^{k_1-1} \sum_{n=0}^{k_2-1} {\sigma_{m,n}}, \; \;
 S_k =  S_{k,k}.\]

\begin {example} \rm
 Let us now describe the basic example of maps satisfying the above conditions.
Let $\{a_i\}_{i=1}^{r_1}$, $\{b_i\}_{i=1}^{r_2}$ be respectively sets of generators
of $F_{r_1}$ and $F_{r_2}$  (the free groups on
$r_1$ and $r_2$ generators; if $r_1 =r_2$ we consider isomorphic copies of the same group)
and let $\{\alpha_i\}_{i=1}^{r_1}$, $\{\beta_i\}_{i=1}^{r_2}$ be sets of $\phi$-invariant
$\star$-automorphisms of the algebra $\alg$, such that $\alpha_j \circ \beta_i = \beta_i
\circ \alpha_j $ for $i\in \{1,\ldots,r_1\}$, $j\in \{1, \ldots, r_2\}$.
Assume that we have group homomorphisms $\Phi_1 : F_{r_1} \longrightarrow Aut(\alg)$ and
$\Phi_2 : F_{r_2} \longrightarrow Aut(\alg)$
defined on the basis elements by $\Phi_1 (a_i) = \alpha_i$, $i \in \{1,\ldots,r_1\}$,
$\Phi_2 (b_j) = \beta_j$, $j \in \{1,\ldots,r_2\}$.
Let for each $n\in \bn$ $w^{(1)}_n$  (respectively $w^{(2)}_n$) denote
a set of reduced words belonging to $F_{r_1}$ ($F_{r_2}$) of length $n$. Further let
$|w^{(1)}_n|$ (respectively $|w^{(2)}_n|$) denote the cardinality of this set (e.g.\ $|w^{(1)}_1|=2r_1$).
The following elements are double-indexed equivalents of objects introduced in \cite{NevoStein}
and will be called the \it Multi Free
Group Actions \rm and the \it Square Free Group Partial Sums:\rm
\[\sigma_{m,n} = \frac{1}{|w^{(1)}_m| |w^{(2)}_n|} \sum_{a \in w^{(1)}_n} \sum_{a \in w^{(2)}_n}
{\Phi_1(a) \circ \Phi_2(b)}, \; \;
S_n = \frac{1}{n^2} \sum_{j=0}^{n-1} \sum_{k=0}^{n-1} {\sigma_{j,k}}. \]
The respective recurrence relations follow from properties
of the free group.
\end{example}

Note that in the situation described in the beginning we can define operators
 $\tilde{\sigma}_{0,1}$
and $\tilde{\sigma}_{1,0}$ acting on  $\Lambda_\phi (\alg_0)$ by
\[ \tilde{\sigma}_{0,1} (\Lambda_\phi (B)) = \Lambda_\phi (\sigma_{0,1} (B)), \; \; B \in \alg_0\]
(similarly $\tilde{\sigma}_{1,0}$).
The complete positivity, unitality and $\phi$-invariance of the initial maps imply that
$\tilde{\sigma}_{0,1}$, $\tilde{\sigma}_{1,0}$
are contractive, and as such can be continuously extended to the whole $H_\phi$ (the extension will
be denoted by the same symbols). Using recurrence relations we define in a natural manner
$\tilde{\sigma}_{0,0}$, $\tilde{\sigma}_{m,n}$ (for $m+n \geq 2$), etc..
    The important fact concerning maps $\tilde{\sigma}_{0,1}$
and $\tilde{\sigma}_{1,0}$ defined in this way is that they are commuting normal operators,
 so they satisfy all  assumptions of theorem
\ref{vNWalker}, except possibly the spectra conditions.

\begin{tw} Let $A \in \alg_0$. If the spectra of $\tilde{\sigma}_{0,1}$ and  $\tilde{\sigma}_{1,0}$
(as operators in $B(H_{\phi})$) are respectively  contained in $D_{p_1}$ and in $D_{p_2}$ then
the multisequence $(S_{m,n} (A))_{m,n=1}^{\infty}$  converges strongly in Pringsheim's sense
to $\hat{A} \in \alg_0$. Moreover if $P \in B(H_{\phi})$ is a projection onto $\{\eta \in H_\phi :
\tilde{\sigma_{0,1}} \eta = \tilde{\sigma_{1,0}}\eta \}$ then $\Lambda_\phi (\hat{A}) = P \Lambda_\phi (A)$.
\label{strong}
\end{tw}

\begin{proof} Theorem \ref{vNWalker} shows that $\tilde{S_{\bt}} \stackrel{\bt\longrightarrow \infty}
    {\longrightarrow} P$ strongly. This implies that for any $\eta_1, \eta_2 \in \alg_0 '$, $A \in \alg_0$
\[ \langle \pi_{\phi} (S_{\bt} (A)) \eta_1 | \eta_2 \rangle = \langle \Lambda_{\phi} (S_{\bt} (A)) | \eta_2
\eta_1^{\flat} \rangle
  = \langle \tilde{S_{\bt}} (\Lambda_{\phi} ( A)) | \eta_2 \eta_1^{\flat} \rangle
\stackrel{\bt\longrightarrow \infty}
    {\longrightarrow}  \langle  P(\Lambda_{\phi} ( A)) | \eta_2 \eta_1^{\flat} \rangle.\]
      If $\psi \in \alg_{\star}$ is defined via $\psi(B) = \langle \pi_{\phi} (B) \eta_1 | \eta_2 \rangle$, $B \in \alg$,
then the absolute value of the right-hand side of the above expression
can be estimated by  $\| \psi\| \cdot \|A\|_{\infty}$ (all $S_{\bt}$ are contractive).
Using the fact that the set
of the above forms is dense in $A_{\star}$ we may conclude that the functional
\[ \alg_{\star} \ni \psi \longrightarrow \lim_{k \longrightarrow \infty} {\psi(S_k (A)} \in \bc\]
is well defined and continuous. Therefore there exists $\hat{A} \in \alg_0$ such that for all
$\psi \in \alg_{\star}$
\[ \psi(S_{\bt}(A)) \stackrel{\bt\longrightarrow \infty} {\longrightarrow} \psi (\hat{A}).\]
It is clear that $\|\hat{A}\|_{\infty} \leq \|A\|_{\infty}$.

Consider again any $\eta_1, \eta_2 \in \alg_0 '$. We have
\[ \langle \pi_{\phi} (\hat{A}) \eta_1 | \eta_2 \rangle = \langle  P \Lambda_{\phi} (A) | \eta_2 \eta_1^{\flat} \rangle
= \langle \pi'_{\phi} (\eta_1)  P \Lambda_{\phi} (A) | \eta_2 \rangle,\]
where $\pi'_\phi$ denotes the right regular representation. As $\eta_2$ was arbitrary, we obtain
\begin{equation} \pi_{\phi} (\hat{A}) \eta_1 = \pi'_{\phi} (\eta_1)  P \Lambda_{\phi} (A). \label{raz}\end{equation}
As all maps $S_{\bt}$ are $*$-preserving,
we can easily prove that $\hat{(A^{\star})} = ({\hat{A}})^{\star}$. This applied to (\ref{raz})
gives
\begin{equation} \pi_{\phi} (\hat{A})^{\star} \eta_1 = \pi'_{\phi} (\eta_1)  P \Lambda_{\phi} (A^{\star}) \label{dwa}.\end{equation}
Using Proposition 10.4 of \cite{Stratila}, we infer from (\ref{raz}) and (\ref{dwa}),
 and the fact that the Hilbert algebra
$\alg_0$ is full, that $P \Lambda_{\phi} (A) \in \alg_0$, $P \Lambda_{\phi} (A) =
\Lambda_{\phi} (\hat{A})$. Now the required strong convergence can be obtained almost immediately,
again taking any $\eta \in \alg_0 '$ :
\[ \pi_{\phi} (S_{\bt} (A)) \eta =  \pi'_{\phi} (\eta) \Lambda_\phi(S_{\bt} (A))
 = \pi'_{\phi} (\eta) \tilde{S_{\bt}} (\Lambda_\phi(A)) \stackrel{\bt\longrightarrow \infty} {\longrightarrow}
     \pi'_{\phi} (\eta)  P \Lambda_{\phi} (A) =  \pi_{\phi} (\hat{A}) \eta \]
and using the fact that $\pi_\phi$ is normal.
\end{proof}

We need the following lemma, which is a straightforward generalisation of Lemma 1
of \cite{NevoStein}.

\begin{lem}
There exist $C_{p_1}, C_{p_2}>0$ such that, for all $m,n \in \bn$ and $A \in \alg^+$,
\[ S_{m,n} (A)\leq C_{p_1} C_{p_2} \frac{1}{3^2 m n} \sum_{j=0}^{3m-1}
 \sum_{l=0}^{3n-1} {\sigma_{0,1}^j \sigma_{1,0}^l (A)}.\]
\label{Nevos} \end{lem}

The following consequence is needed below.

\begin{lem} For any $A \in \alg$ and $\bm \in {\bn}^2$ the multisequence
$\| S_{\bt} (S_{\bm} (A) - A )\|_{\infty}$ tends to $0$ as $\max {\bt}$ tends to $\infty$ (and so
also in Pringsheim's sense).  \label{comb}
\end{lem}

\begin{proof} We begin with the following observation:
\[ S_{\bm}(A) - A = \sum_{j=0}^{m_1-1} \sum_{l=0}^{m_2-1}
 {\left(\frac{1}{m_1 m_2}
 (\sigma_{j,l}(A) - A) \right)},\]
 so it is enough to prove convergence for expressions such as
$\| S_{\bt} (\sigma_{j,l} (A) - A )\|_{\infty}$. In turn we can reduce this to proving that
for each $j, l\in  \bn$
$\| S_{\bt} (\sigma_{0,1}^j \circ \sigma_{1,0}^l(A) - A )\|_{\infty}$ tends to $0$ as $\max{\bt}$ tends to $\infty$.
 However this can be obtained with the help of the previous lemma by considering standard Cesaro
averages.
\end{proof}

We will use the existence of a convenient
decomposition of a selfadjoint operator in $\alg_0$, proved by D.Petz in
\cite{Petz1}:

\begin{lem} Suppose that $B \in \alg_0$, $B=B^{\star}$. Then there exist $C\in \alg$, $C=C^{\star}$,
$D,E \in \alg^+$ such that $B=C+D-E$, $\|C\|_{\infty} \leq \phi(B^2)^{\frac{1}{2}}$, $\phi(D) \leq
\phi (B^2)^{\frac{1}{2}}$, $\phi(E) \leq \phi (B^2)^{\frac{1}{2}}$ and $\|C\|_{\infty}, \|D\|_{\infty},
\|E\|_{\infty} \leq \|B\|_{\infty}$. \label{decom}
\end{lem}

Now we can formulate the first of the two main results of this section.

\begin{tw} Let $A \in \alg_0$. If the spectra of $\tilde{\sigma}_{0,1}$ and  $\tilde{\sigma}_{1,0}$
(as operators in $B(H_{\phi})$) are respectively  contained in $D_{p_1}$ and in $D_{p_2}$ then
the sequence $(S_n (A))_{n=1}^{\infty}$ is b.a.u.\
convergent to $\hat{A} \in \alg_0$. \label{bauc}
\end{tw}

\begin{proof} Theorem \ref{strong} implies that if $k\in \bn$ and
\[ \xi_k = \tilde{S_k}(\Lambda_\phi(A)) - P \Lambda_\phi(A) \]
then $\|\xi_k\|_2$ tends to $0$ as $k$ tends to $\infty$.
As $\hat{A}$ is invariant under $\sigma_{1,0}$ and $\sigma_{0,1}$,
we have $S_k(\hat{A}) = \hat{A}$ for all $k \in \bn$.
Moreover $\xi_k \in \Lambda_\phi(\alg_0)$
and if $\xi_k = \Lambda_\phi(B_k)$, $B_k \in \alg_0$, we have
\[A - \hat{A}  = B_k + A - S_k(A), \; \; \phi(B_k^{\star} B_k) \stackrel{k\longrightarrow \infty} {\longrightarrow} 0. \]
Decomposing $A$ into its real and imaginary part we can  assume that $B_k = B_{k}^{\star}$.
Let us fix $\epsilon>0$ and choose a subsequence $(k_n)_{n=1}^{\infty}$
such that $\phi(B_{k_n}^2)^{\frac{1}{2}} \leq  n^{-1} 2^{-n-1} \epsilon$.
For each $n\in \bn$ we can decompose $B_{k_n}$ according
to Lemma \ref{decom}, $B_{k_n} = C_{k_n}+D_{k_n} - E_{k_n}$.
Without loss of generality we assume that say $E_{k_n} = 0$.
Now we apply Lemma \ref{Petzmax} (or rather actually its version
in Remark \ref{comment})  for maps $\sigma_{0,1}$, $\sigma_{1,0}$ and a
sequence $(D_{k_n})_{n=1}^{\infty}$, with estimation numbers respectively
equal to $\frac{1}{n}$, as a result finding a projection $p \in P_\alg$ such that
\[ \phi(p^{\perp}) \leq  2 \sum_{n=1}^{\infty }{n n^{-1} 2^{-n-1} \epsilon} = \epsilon, \]
\[ \|p \left(\frac{1}{r^2}\sum_{l_1=0}^{r-1} \sum_{l_2 =0}^{r-1}
{\sigma_{1,0} ^{l_1} \circ \sigma_{0,1}^{l_2} (D_{k_n})} \right) p \|_{\infty} \leq 2 n^{-1}, \; \; r, n \in \bn,\]

In the end, using Lemma \ref{Nevos} we obtain (for any $n, k \in \bn$)
\[ \| p (S_k(A - \hat{A})) p\|_{\infty} = \| p (S_{k} ( B_{k_n} + A - S_{k_n}(A)))p\|_{\infty} \leq \]
\[ \| p (S_{k} B_{k_n})  p \|_{\infty} + \| p (S_{k} (A - S_{k_n}(A))p\|_{\infty} \leq \]
\[ \|B_{k_n}\|_2 + C_{p_1} C_{p_2} \chi_2 \frac{2}{n} +  \| S_{k} (A - S_{k_n}(A))\|_{\infty},\]
and an application of Lemma \ref{comb} ends the proof. \end{proof}

The scheme described in section 2 allows us to deduce immediately the
second important result.

\begin{tw} Let $M$ be a von Neumann algebra with a normal semifinite faithful
trace $\tau$ and let  $x \in L^1(M)$. If $(\sigma_{\bt})_{\bt \in \bn^2}$ is the sequence of maps
acting on $M$ and satisfying the conditions described before Theorem \ref{strong}
then the sequence $(S_n(x))_{n=1}^{\infty}$ converges b.a.u.\ to some
 $\widehat{x} \in L^1(M)$.
\end{tw}

\begin{proof}  Assume that $x\geq 0$.
As it is clear that $\sigma_{0,1} $ and $\sigma_{1,0}$ are commuting kernels,
we can (as was done above) use Lemma \ref{Petzmax} and Lemma \ref{Nevos}
to deduce that for every
$\epsilon >0$ there exists
$p \in P_M$ such that $\tau(p^{\perp}) < \epsilon $ and for all $n \in \bn$
\[ \| p S_n (x) p \|_{\infty} \leq  \epsilon^{-1} C_{p_1} C_{p_2}\chi_2 \|x\|_1.\]
Obviously each $S_n$ treated as a map from $L^1(M)_{\rm sa}$ to
$\widetilde{M}$ is positive and
continuous. Theorem \ref{bauc} implies that for any $x \in M^+\cap L^2(M)$ the sequence
$(S_n(x))_{n=1}^{\infty}$ is b.a.u. convergent.
 As $M^+\cap L^2(M)$ is a minorantly dense
subset of $L^1(M)_{\rm sa}$, we are in position to apply the noncommutative Banach
principle (Theorem \ref{Banach}) to end the proof. \end{proof}

As a special case, putting $p_1 = \ldots = p_d =1$ we obtain the
 noncommutative generalization of the classical
result of A.Brunel:

\begin{cor}
Assume that $\alpha_1, \ldots, \alpha_d$ are commuting, normal,
completely positive, unital, $\tau$-invariant maps acting on $M$.
Then for each $x \in L^1(M)$ the sequence $(s_n(x))_{n=1}^{\infty}$,
\[s_n (x) = \frac{1}{n^d} \sum_{i_1=0}^{n=1} \ldots \sum_{i_d=0}^{n-1}
\alpha_1^{i_1} \circ \ldots \alpha_d^{i_d} (x), \; \; n\in \bn\]
is b.a.u.\ convergent.
\end{cor}

All the results remain  true if instead of considering the averages over
squares we deal with so-called sequences of indices tending to
infinity but \it remaining
in a sector of $\bn^d$.\rm This means that we consider averaging over sets
of the type $\{1,\ldots, k_1(n)\} \times \ldots \times \{1,\ldots, k_d(n)\}$,
for which there exists $C>0$ such that $\frac{k_i(n)}{k_j(n)} <C$
for all $i,j \in \{1,\ldots,d\}$, $n\in \bn$.

\vspace{1 cm}
\noindent
{\bf ACKNOWLEDGEMENTS}

During the preparation of this work the author was partially supported by the
KBN Research Grant 2P03A 030 24 and by the European Comission HPRN-CT-2002-00279,
RTN QP-Applications.

The author would like also to express his gratitude to Vladimir Chilin and
Semyon Litvinov, whose remarks essentially improved the final form of this
paper.

\end{document}